 \documentclass[final,1p,times,authoryear]{elsarticle}

\usepackage{epsfig}

\usepackage{amssymb}
\usepackage{amsthm}
\usepackage{amsmath}
\usepackage{amsfonts}
\usepackage{mathtools}
\usepackage{booktabs,lipsum}
\usepackage[colorlinks=true,linkcolor=red]{hyperref}

\usepackage{lineno}
\biboptions{sort&compress}

\journal{Int. J. Appl. Comput. Math}

\begin{document}

\begin{frontmatter}



\title{Numerical Study of Astrophysics Equations by Meshless Collocation Method Based on Compactly Supported Radial Basis Function}


\author[mymainaddress]{K. Parand\corref{mycorrespondingauthor}}
\cortext[mycorrespondingauthor]{Corresponding author. Member of research group of Scientific Computing. Fax:
+98 2122431653.}
\ead{k\_parand@sbu.ac.ir}

\author[mymainaddress]{M. Hemami}
\ead{mohammadhemami@yahoo.com}

\address[mymainaddress]{Department of Computer Sciences, Shahid Beheshti University, G.C. Tehran 19697-64166, Iran}

\begin{abstract}

In this paper, we propose compactly supported radial basis  functions for solving some well-known classes of
astrophysics problems categorized as non-linear singular initial ordinary differential equations
on a semi-infinite domain. To increase the convergence rate and to decrease the collocation
points, we use the compactly supported radial basis function through the integral operations. Afterwards,
some special cases of the equation are presented as test examples to show the reliability of the
method. Then we compare the results of this work with some results and show that the
new method is efficient and applicable.

\end{abstract}


\end{frontmatter}

 \section{Introduction}
\label{intro}
In recent decades, the so-called meshless methods have been extensively used to find approximate solutions of various types of linear and non-linear equations \citep{MeshfreeFasshhauer} such as differential equations (DEs) and integral equations (IEs). Unlike the other methods which were used to mesh the domain of the problem, meshless  methods don't require a structured grid and only make use of a scattered set of collocation points regardless of the connectivity information between the collocation points.\par
For the last years, the radial basis functions (RBFs) method was known as a powerful tool for the scattered data interpolation problem.
One of the domain-type meshless methods is given in \citep{Kansa} in 1990, which directly collocates RBF, particularly the multiquadric (MQ) to find an approximate solution of linear and non-linear DEs. The RBFs can be compactly and globally supported. The global RBF's are infinitely differentiable and contain a free parameter $c$, called the shape parameter \citep{Islam,Buhman2,Sara}. The interested reader is referred to the recent books and paper by Buhmann \citep{Buhman2,Buhman3} and Wendland \citep{Wendland2} for more basic details about RBFs, compactly and globally supported and the convergence rate of the RBFs. In RBF method  using globally RBFs, if $c$ increases, the system of equations to be solved becomes  ill-conditioned. Cheng et al. \citep{Cheng} showed that when $c$ is very large then the RBFs system error is of exponential convergence. To overcome the problems of globally RBF, alternative ways suggested such as domain decomposition\citep{Li}, LU decomposition\citep{KansaHon}, good matrix pre-conditioners\citep{KansaHon} change the global support RBFs with compactly support RBFs (CSRBFs)\citep{Shen} which are local support and have not any parameter. \par
There are two basic approaches for obtaining basis functions from RBFs, namely direct approach (DRBF) based on differential process \citep{Kansa2} and indirect approach (IRBF) based on an integration process \citep{Duy1,Duy2,Duy3}. Both approaches were tested on the solution of second order DEs and the indirect approach was found to be superior to the direct approach \citep{Duy1}.\par
In this paper, we use the indirect CSRBF (ICSRBF) for finding the solution of Lane-Emden type equations and also Isothermal gas sphere, White  Dwarf equation.\par
This paper is arranged as following:\\
 In Section 2 we describe the Lane-Emden equations. In Section 3  we survey several methods that have been used to solve Lane-Emden type equations. In Section 4, the properties of CSRBF and the way to construct the  ICSRBF method for this type of equations are described. In Section 5 the proposed method is applied to some types of Lane-Emden equations, and a comparison is made with the existing solutions that were reported in other published works. Finally we give a brief conclusion in the last section.
 \section{Lane-Emden type equations}
\label{sec:1}
The Lane-Emden equation describes a variety of phenomena in theoretical physics and astrophysics, including aspects of stellar structure, the thermal history of a spherical cloud of gas, isothermal gas spheres, and thermionic currents.\\\par
Let $P(r)$ denote the total pressure at a distance $r$ from the center of spherical gas cloud. The total pressure is due to the usual gas pressure and a contribution from radiation:
\begin{equation}
P=\frac{1}{3}\varsigma T^4+\frac{RT}{v},\nonumber
\end{equation}
where $\varsigma$ , $T$, $R$  and $v$ are respectively the radiation constant, the absolute temperature, the gas constant, and the specific volume\citep{Ravi}. Let $M(r)$ be the mass within a sphere of radius $r$ and $G$ the constant of gravitation. The equilibrium equations for the configuration are
\begin{equation}
\label{dens}
\frac{dP}{dr}=-\frac{GM(r)}{r^2}, \\
\end{equation}
\begin{equation}
\frac{dM(r)}{dr}=4\pi \rho r^2,\nonumber
\end{equation}
where $\rho$ is the density, at a distance $r$ from the center of a spherical star.\\
Eliminating $M$ yields:
\begin{equation}
\frac{1}{r^2}\frac{d}{dr}\Bigg(\frac{r^2}{\rho}\frac{dP}{dr}\Bigg)=-4\pi G\rho .\nonumber
\end{equation}

Pressure and density $\rho=v^{-1}$ vary with $r$ and $P=K\rho^{1+\frac{1}{m}} $ where $K$ and $m$ are constants.\\
We can insert this relation into Eq. (\ref{dens}) for the hydrostatic equilibrium condition and from this rewrite equation to:
\begin{equation}
\Bigg[\frac{K(m+1)}{4\pi G}\lambda ^{\frac{1}{m}-1}\Bigg]\frac{1}{r^2}\frac{d}{dr}\Bigg(r^2\frac{dy}{dr}\Bigg)=-y^m,\nonumber
\end{equation}
where $\lambda$ representing the central density of the star and $y$ is dimensionless quantity that are both related to $\rho$ through the following relation:
\begin{equation}
\rho=\lambda y^m,\nonumber
\end{equation}
and let:
\begin{eqnarray}
&&r=ax,\nonumber\\
&&a=\Bigg[\frac{K(m+1)}{4\pi G}\lambda^{\frac{1}{m}-1}\Bigg]^{\frac{1}{2}}.\nonumber
\end{eqnarray}
Inserting these relations into our previous relations we obtain the Lane-Emden equation \citep{Dehghan1}:
\begin{equation}
\frac{1}{x^2}\frac{d}{dx}\Bigg(x^2\frac{dy}{dx}\Bigg)=-y^m,\nonumber
\end{equation}
and with simplifying previous equation we have:
\begin{equation}
y''(x)+\frac{2}{x}y'(x)+y^m(x)=0,~~~x>0,\nonumber
\end{equation}
with the boundary conditions:
\begin{equation}
y(0)=1,~~~y'(0)=0.\nonumber
\end{equation}
It has been claimed in the literature that only for $m=0,~1$ and $5$ the solutions of the Lane-Emden equation could be exact. For the other values of $m$, the Lane-Emden equation is to be integrated numerically.
 In this paper, we solve it for $m=1,1.5,2,2.5,3,4$ and $5$.\\
 \section{Methods have been proposed to solve Lane-Emden type equation}
\label{sec:3}
in years, many analytical and numerical methods have been used to solve Lane-Emden equations, the main difficulty arises in the singularity of the equations at $x=0$. Currently, most techniques which were used in handling the Lane-Emden type problems are based on either series solutions or perturbation techniques.\\
Bender et al. \citep{Bender}, proposed a new perturbation technique based on an artificial parameter $\delta$, the method is often called $\delta-$method.\par
Mendelzweig and Tabakin \citep{Mendelzweig} used quasi-linearization approach to solve the standard Lane-Emden equation. This method approximates the solution of a non-linear differential equation by treating the non-linear terms as a perturbation about the linear ones, and unlike perturbation theories is not based on the existence of some small parameter. He showed that the quasi-linearization method gives excellent results when applied to different no-nlinear ordinary differential equations in physics, such as the Blasius, Duffing, Lane-Emden and Thomas-Fermi equations.\par
Shawagfeh \citep{Shawagfeh} applied a non-perturbative approximate analytical solution for the Lane-Emden type equation using the Adomian decomposition method. His solution was in the form of a power series. He used Pad\'e approximants method to accelerate the convergence of the power series.\par
In \citep{Wazwaz1}, Wazwaz employed the Adomian decomposition method with an alternate framework designed to overcome the difficulty of the singular point. It was applied to the differential equations of Lane-Emden type. In future \citep{Wazwaz2} he used the modified decomposition method for solving analytical treatment of non-linear differential equations such as Lane-Emden equations. The modified method accelerates the rapid convergence of the series solutions, dramatically reduces the size of work, and provides the solution by using few iterations only without any need to the so-called Adomian polynomials.\par
Liao \citep{LiaO} provided a reliable, easy-to-use analytical algorithm for Lane-Emden type equations. This algorithm logically contained the well-known Adomian decompositions method. Different from all other analytical techniques, this algorithm itself provides us with a convenient way to adjust convergence regions even without Pad\'e technique.\par
He \citep{He} employed Ritz's method to obtain an analytical solution of the problem. By the semi-inverse method, a variational principle is obtained for the Lane-Emden equation, which he gave much numerical convenience when applied to finite element methods or the Ritz method.\par
Parand et al. \citep{Parand1, Parand2,Parand15} presented some numerical techniques to solve higher ordinary differential equations such as Lane-Emden. Their approach was based on a rational Chebyshev and rational Legendre tau method. They presented the derivative and product operational matrices of rational Chebyshev and rational Legendre functions.\par
These matrices together with the tau method were utilized to reduce the solution of these physical problems to the solutions of systems of algebraic equations. Also, in some paper's \citep{Parand7,Parand8,Parand9,Parand12}, Parand et al. applied Hermite function collocation method (HFC) , Bessel function collocation method and meshless collocation method based on Radial basis function (RBFs) to solving the Lane-Emden type equations.\par
Ramos \citep{Ramos1,Ramos2,Ramos3,Ramos4} solved Lane-Emden equation through different methods. In \citep{Ramos1} he presented linearzation methods for singular initial value problems in second order ordinary differential equation such as Lane-Emden. These methods result in linear constant-coefficient ordinary differential equations which can be integrated analytically, thus they yield piecewise analytical solutions and globally smooth solutions \citep{Ramos2}. Later, he developed piecewise-adaptive decomposition methods for the solution of non-linear ordinary differential equations. Piecewise-decomposition methods provide series solutions in intervals which are subject to continuity conditions at the end points of each interval and their adoption is based on the use of either a fixed number of approximants and a variable step size, a variable number of approximants and a fixed step size or a variable number of approximants and a variable step size. In \citep{Ramos3}, series solutions of the Lane-Emden equation, based on either a volterra integral equation formulation or the expansion of the dependent variable in the original ordinary differential equation are presented and compared with series solutions obtained by means of integral or differential equations based on a transformation of the dependent variables.\par
Yousefi \citep{Yousefi} used integral operator and converted Lane-Emden equations to integral equations and then applied Legendre Wavelet approximations. In this work properties of Legendre wavelet together with the Gaussian integration method were utilized to reduce the integral equations to the solution of algebraic equations. By his method, the equation was formulated on [0, 1].\par
Chowdhury and Hashim \citep{Hashim} obtained analytical solutions of the generalized Emden-Fowler type equations in the second order ordinary differential equations by homotopy-perturbation method (HPM). This method is a coupling of the perturbation method and the homotopy method. The main feature of the HPM \citep{Dehghan2} is that it deforms a difficult problem into a set of problems which are easier to solve. HPM yields solution in convergent series forms with easily computable terms.\par
Aslanov \citep{Aslanov} constructed a recurrence relation for the components of the approximate solution and investigated the convergence conditions for the Emden-Fowler type of equations. He improved the previous results on the convergence radius of the series solution.\par
Dehghan and Shakeri \citep{Dehghan3} investigated Lane-Emden equation using the variational iteration method 
and showed the efficiency and applicability of their procedure for solving this equation. Their technique does not require any discretization, linearization or small perturbations and therefore reduces the volume of computations.\par
Bataineh et al. \citep{Bataineh} obtained analytical solutions of singular initial value problems (IVPs) of the Emden-Fowler type by the homotopy analysis method (HAM). Their solutions contained an auxiliary parameter which provided a convenient way of controlling the convergence region of the series solutions. It was shown that the solutions obtained by the Adomian decomposition method (ADM) and the homotopy-perturbation method (HPM) are only special cases of the HAM solutions.\par
 singh et al. \citep{Singh} used the modified Homotopy analysis method for solving the Lane-Emden-equation and White Dwarf equation.\par
Adibi and Rismani in \citep{Adibi} proposed the approximate solutions of singular IVPs of the Lane-Emden type in second-order ordinary differential equations by improved Legendre-spectral method. The Legendre-Gauss point used as collocation nodes and Lagrange interpolation is employed in the Volterra term.\par
Karimi vanani and Aminataei \citep{Vanani}  provide a numerical method which produces an approximate polynomial solution for solving Lane-Emden equations as singular initial value problem. They are first, used as an integral operator and convert Lane-Emden equations into integral equations, then convert the acquired integral equation into a power series and finally, they transforming the power series into pade series form.\par
Kaur et al. \citep{Kaur}, obtained the Haar wavelet approximate solutions for the generalized Lane-Emden equations. This method was based on the quasi-linearization approximation and replacement of an unknown function through a truncated series of Haar wavelet series of the function.\par
Other researchers try to solve the Lane-Emden type equations with several methods, For example, Y{\i}ld{\i}r{\i}m and \"{O}zi\c{s} \citep{Ozis1,Ozis2} by using HPM and VIM methods, Benko et al. \citep{Benko} by using Nystr\"{o}m method, Iqbal and javad \citep{IQbal} by using Optimal HAM, Boubaker and Van Gorder \citep{Boubaker} by using Boubaker polynomials expansion scheme, Da\c{s}c{\i}o\u{g}lu and Yaslan \citep{Dascioglu} by using Chebyshev collocation method, Y\"{u}zba\c{s}{\i} \citep{Yuzbasi1,Yuzbasi2} by using Bessel matrix and improved Bessel collocation method, Boyd \citep{Boyd} by using Chebyshev spectral method, Bharwy and Alofi \citep{Alofi} by using Jacobi-Gauss collocation method, Pandey et al. \citep{Pandey1,Pandey2} by using Legendre and Brenstein operation matrix, Rismani and monfared \citep{Rismani} by using Modified Legendre spectral method, Nazari-Golshan et al. \citep{Golshan} by using Homotopy perturbation with Fourier transform, Doha et al. \citep{Doha} by using second kind Chebyshev operation matrix algorithm, Carunto and bota \citep{bota} by using Squared reminder minimization method, Mall and Chakaraverty \citep{Mall} by using Chebyshev Neural Network based model, G\"{u}rb\"{u}z and sezer \citep{Gurbuz} by using Laguerre polynomial, Kazemi-Nasab et al. \citep{Kazemi} by using Chebyshev wavelet finite difference method, Hosseini and Abbasbandy \citep{Abbas1} by using combination of spectral method and ADM method and Azarnavid et al. \citep{Abbas2} by using Picard-Reproducing Kernel Hilbert Space Method\par
\section{ ICSRBF method}
\label{sec:4}
\subsection{CSRBF}
\label{ssec:41}
Many problems in science and engineering arise in infinite and semi-infinite domains. Different numerical methods have been proposed for solving problems on various domains such as FEM\citep{DH10,DH11}, FDM\citep{DH9,DH10}, Spectral\citep{DH5,Parand8,Parand13,Parand14,Parand17} methods and meshfree methods\citep{DH3,DH4,dehghan2009meshless,Abbas3}.
The use of the RBF is one of the popular meshfree method for solving the differential equations  \citep{DH1,DH2,radi1,radi2}. For many years the global radial basis functions such as Gaussian, Multi quadric, Thin plate spline, Inverse multiqudric and  others were used  \citep{DH6,DH7,DH8} to solve different DEs and interpolation. These functions are globally supported and generate a system of equations  with ill-condition full matrix. To convert the ill-condition matrix to a well-condition matrix, CSRBFs can be used  instead of global RBFs. CSRBFs can convert the global scheme into a local one with banded matrices, Which makes the RBF method more feasible for solving large-scale problem \citep{Wong}.
\subsection*{Wendland's functions}
\label{ssec:42}
The most popular family of CSRBF are Wendland functions. These functions were introduced by Holger Wendland in 1995 \citep{wendland}. he started with the truncated power function $\phi_l(r)=(1-r)^l_+$ which is strictly positive definite and radial on $\mathbb{R}^s$ for $l\geq \lfloor\frac{s}{2}\rfloor+1$ , and then he walks through dimensions by repeatedly applying the mont\'{e} operator (I).\\
\textbf{Definition 1 \citep{MeshfreeFasshhauer}}~~
with $\phi_l(r)=(1-r)^l_+$ he defines
 \begin{equation}
 \phi_{s,k}=I^k\phi_{\lfloor\frac{s}{2}\rfloor +k+1},
 \end{equation}
 it turns out that the functions $\phi_{s,k}$ are all supported on [0,1].\\
\textbf{Theorem 1 \citep{MeshfreeFasshhauer}}~~
 The functions $\phi_{s,k}$ are strictly positive definite (SPD)  and radial on $\mathbb{R}^s$ and are of the form
 \begin{equation}
\phi _{s,k}(r)=\begin{cases}
p_{s,k}(r)&   r\in [0,1],\nonumber\\
0&       r>1,
\end{cases}
\end{equation}
with a univariate polynomial $p_{s,k}$ of degree $\lfloor \frac{s}{2}\rfloor+3k+1$. Moreover, ،$\phi_{s,k}\in C^{2k}(R)$ are unique up to a constant factor, and the polynomial degree is minimal for given space dimension $s$ and smoothness $2k$ \citep{MeshfreeFasshhauer}.
Wendland gave recursive formulas for the functions $\phi_{s,k}$ for all $s, k$.\\
\textbf{Theorem 2 \citep{MeshfreeFasshhauer}}
~~~~The functions  $\phi_{s,k}$, $k=0,~1,~ 2,~ 3,$ have form
 \begin{eqnarray}
&&\phi_{s,0}=(1-r)^l_+,\nonumber\\
&&\phi_{s,1}\doteq(1-r) _+^{l+1}[(l+1)r+1],\nonumber\\
&&\phi_{s,2}\doteq(1-r)_+^{l+2}[(l^2+4l+3)r^2+(3l+6)r+3],\nonumber\\
&&\phi_{s,3}\doteq(1-r)_+^{l+3}[(l^3+9l^2+23l+15)r^3+(6l^2+36l+45)r^2\nonumber\\
&&+(15l+45)r+15],\nonumber
\end{eqnarray}
where $l=\lfloor \frac{s}{2}\rfloor +k+1$, and the symbol $\doteq$ denotes equality up to a multiplicative positive constant.\\
~~~~~~The case $k=0$ directly follows  from the definition 1. application of the definition 1 for the case $k=1$ yields
\begin{eqnarray}
&&\phi_{s,1}(r)=(I\phi_l)(r)\nonumber
=\int^\infty_r t\phi_l(t)dt  \nonumber\\
&&=\int^\infty _r t(1-t)^l_+dt\nonumber
=\int^1_r t(1-t)^ldt\nonumber\\
&&=\frac{1}{(l+1)(l+2)} (1-r)^{l+1}[(l+1)r+1],\nonumber
\end{eqnarray}
where the compact support of $\phi_l$ reduces the improper integral to a definite integral which can be evaluated using integration by parts. The other two cases are obtained similarly by repeated application of $I$.\citep{MeshfreeFasshhauer}
We showed most of the wendland functions in table \ref{Table. 1.} .
\begin{table}[htbp]
\caption{Wendland's compactly supported radial function for various choices of k and s=3.}
\label{Table. 1.}       
\centering  \begin{tabular}{lll}
\hline\noalign{\smallskip}
 $\phi_{s,k}$ & smoothness & SPD   \\
\noalign{\smallskip}\hline\noalign{\smallskip}
 $\phi_{3,0}(r) =(1-r)^2_+$ & $C^0$ & $\mathbb{R}^3$  \\
 $\phi_{3,1}(r)\doteq(1-r)^4_+(4r+1)$ & $C^2$ & $\mathbb{R}^3$ \\
 $\phi_{3,2}(r)\doteq(1-r)^6_+(35r^2+18r+3)$ & $C^4$ & $\mathbb{R}^3$  \\
 $\phi_{3,3}(r)\doteq(1-r)^8_+(32r^3+25r^2+8r+1)$ & $C^6$ & $\mathbb{R}^3$  \\
 $\phi_{3,4}(r)\doteq(1-r)^{10}_+(429r^4+450r^3+210r^2+50r+5)$ & $C^8$ & $\mathbb{R}^3$ \\
 $\phi_{3,5}(r)\doteq(1-r)^{12}_+(2048r^5+2697r^4+1644r^3+566r^2+108r+9)$ & $C^{10}$ & $\mathbb{R}^3$  \\
\noalign{\smallskip}\hline
\end{tabular}
\end{table}
\subsection{Interpolation by CSRBFs}
\label{ssec:43}
To interpolate or approximate one dimentional function $y(x)$, we can represent it by a CSRBF as

\begin{equation}
y(x)\approx y_n(x)=\sum_{i=1}^N \xi_i \phi_i (x)=\Phi ^T(x)\Xi\nonumber,
\label{systmat}
\end{equation}
where
\begin{eqnarray}
&&\phi_i (x)=\phi (\dfrac{\|x-x_i\|}{r_\omega}),\nonumber\\
&&\Phi ^T(x)=[\phi_1 (x),\phi_2 (x),\cdots ,\phi _N(x)],\nonumber\\
&&\Xi =[\xi_1,\xi_2, \cdots , \xi_N ]^T,\nonumber
\end{eqnarray}
\begin{equation}
\end{equation}
$x$ is the input, $r_\omega$ is the local support domain and $\xi_i$s are the set of coefficients to be determined. By using the local support domain, we mapped the domain of problem to CSRBF local domain. By choosing $N$ interpolate points $(x_j,~ j=1, 2,\cdots,  N)$ in domain:
\begin{equation}
y_j=\sum_{i=1}^N\xi_i\phi_i (x_j)   (j=1, 2, \cdots, N).\nonumber
\end{equation}
To summarize the discussion on the coefficients matrix, we define
\begin{equation}
A\Xi = Y,
\label{system}
\end{equation}
where :
\begin{eqnarray}
&&Y=[y_1, y_2, \cdots,  y_N]^T,\nonumber\\
&&A=[\Phi ^T(x_1), \Phi ^T(x_2), \cdots , \Phi ^T(x_N)]^T\nonumber\\
&&=\begin{pmatrix} 
\phi_1(x_1)&\phi_2(x_1)&\cdots &\phi_N(x_1)\cr \phi_1(x_2)&\phi_2(x_2)&\cdots &\phi_N(x_2)\cr \vdots &\vdots &\ddots &\vdots \cr \phi_1(x_N)&\phi_2(x_N)&\cdots &\phi_N(x_N)\nonumber
\end{pmatrix}.
\end{eqnarray}
Note that $\phi_i(x_j)=\phi(\dfrac{\|x_i-x_j\|}{r_\omega})$, by solving the system (\ref{system}), the unknown coefficients $\xi_i$  will be achieved.
\subsection{ICSRBF method}
\label{ssec:44}
In this paper, we construct the $\phi_i(x)$ in Eq. (\ref{Original}) by using the Wendland function with parameter $k=3$ and $l=5$.\\
At first, we approximate the highest order derivative in the problem by expansion of CSRBFs:\\
\begin{equation}
\dfrac{d^2}{dx^2}y(x)\simeq y^{(2)}_n(x)=\sum_{i=1}^N \xi_i \phi_i(x), 
\label{Original}
\end{equation}
where $\phi_i (x)$ is the CSRBF and $r_\omega$ is the local support domain and  coefficients $\xi_i$ must be determined. Now we define the lower order of derivatives and unknown function $y(x)$ using Gauss-Legendre integration as follows:
\begin{eqnarray}
&&\dfrac{d}{dx}y(x)\simeq\int^x_0 y^{(2)}_n(t)dt-\dfrac{d}{dx}y(x)|_{x=0}=\dfrac{x}{2}\sum_{j=1}^q\sum_{i=1}^N w[j]\xi_i \phi_i (\dfrac{x}{2}\eta[j]+\dfrac{x}{2})=(y^{(1)}_n(x)),\\
&&y(x)\simeq   \int_0^x (y^{(1)}_n(t))dt-y(0)=(y_n(x)),
\end{eqnarray}
where $w[j]$ are weighted coefficients of Gauss-Legendre integration and define as follow:
\begin{equation}
w[j]=\frac{2}{(1-\eta[j]^2)(\frac{d^2}{dx^2}P_q(x))|_{x=\eta[j]}},\nonumber
\end{equation}
$P_q(x)$ is Legendre polynomial of order $q$ and $\eta[j]$ is $j$-th root of $P_{q-1}(x)$.\\
To satisfy the initial conditions of problem, we have:
\begin{eqnarray}
&&(y^{(1)}_n(x))=  \int_0^x y^{(2)}_n(t)dt,\\
&&(y_n(x))=  \int_0^x (y^{(1)}_n(t))dt+1,
\end{eqnarray}
then, we produce the residual function as follows:
\begin{equation}
Res(x)=xy^{(2)}_n(x)+\alpha (y^{(1)}_n(x))+x(y_n(x))^m,
\end{equation}
now, obtain $N$-nodes as follows:
\begin{equation}
\label{formula}
x_j=L(\frac{i}{N})^\gamma ,~~~j=1, 2, \cdots,  N,
\end{equation}
where $L$ is the last collocate node in domain and $\gamma$ is arbitrary parameter. For  Lane-Emden type equations, we selected $\gamma$ between 1.5 and 1.8, because with these values we found the best $\|Res\|_2$ for solution equations.\\
By using ${x_j},  j=1, 2, ..., N $ and residual function, we obtain $N$ equation, so by solving these equations we obtain the unknown coefficients ${\xi_i}_{i=1}^N$.\\
The result of this section can be summarized in the following algorithm for the IVP:
\begin{equation}
F(x,y(x),y'(x),y''(x))=0,~~~~y(0)=a,~~~y'(0)=b.
\end{equation}
\textbf{Algorithm} The algorithm works in the following manner:\\
(1) Choose $N$ center points $\{X_j\}^N_{j=1}$ from domain $[0, \infty)$.\\
(2) Approximate $y''(x)$ as the from $u''_N(x)=\sum_{i=1}^N \xi_i\phi_i(x)$.\\
(3) Obtain $y'(x)$ by using defined integral operation $I_\chi(h(x))=\int_0^x h(t)dt$ in the form $u'_N(x)= \sum_{i=1}^N \xi_i \int _0^x \phi_i(t) dt+b$.\\
(4) Obtain $y(x)$ by using defined integral operation $I_\chi(h(x))=\int_0^x h(t)dt$ in the form $u_N(x)=  \int _0^x u'_N(t) dt+a$.\\
(4) Substitiute $u_N(x)$, $u'_N(x)$ and $u''_N(x)$ into the main problem and create residual function Res(x).\\
(5) Substitiute collocation points $\{X_j\}_{j=1}^{N}$ into the Res(x).\\
(6) Solve the $N$ equations with $N$ unknown coefficients $\{\xi_i\}_{i=1}^N$ and  find the numerical solution.

\begin{equation}
Res(x_j)=0,~~~j=1, 2, \cdots,  N.\nonumber
\end{equation}

\section{Application}
\label{sec:5}
In this section, we apply ICSRBF method to solve the Lane-Emden type equations. In general the Lane-Emden type equation are formulated as follows:
\begin{equation}
y''(x)+\frac{\alpha}{x}y'(x)+p(x)q(y(x))=h(x),~~~\alpha x\geq 0,
\label{lane-eq}
\end{equation}
with initial conditions :
\begin{equation}
y(0)=A ,~~~y'(0)=B,
\end{equation}\\

where $\alpha$, $A$ and $B$ are real constants and $p(x)$, $q(y(x))$ and $h(x)$ are some given function. To apply the collocation method, we construct the residual function by substituting $y_n(x)$  in the Lane-Emden type Eq. (\ref{lane-eq}):
\begin{equation}
Res(x)=y^{(2)}_n(x)+\frac{\alpha}{x}y^{(2)}_n(x)+p(x)q(y_n(x))-h(x)\nonumber.
\end{equation}
The equation for obtaining the coefficient $\xi_is$ arise from equalizing $Res(x)$ to zero at $N$ point (\ref{formula}):
\begin{equation}
Res(x_j)=0,~~~j=0,1,2,...,N.
\end{equation}
By solving this set of equations we obtain a approximate function $y_n(x)$. Note that these $N$ equations generate a set of $N$ nonlinear equations which can be solved by a well-known method such as Newton method for unknown coefficient $\xi_is$. 
\subsection{Example 1 (The standard Lane-Emden equation)}
\label{sssec:51}
\begin{table}[htbp]
\caption{Comparison of the first zeros of standard Lane Emden equations between the present method and exact numerical value given by Horedt\citep{Horedt}.}
\label{Table. 2.}   
\centering   \begin{tabular}{lllllll}
\hline\noalign{\smallskip}
 m &  N & $r_\omega$ & L & The present method & Horedt & Error   \\
\noalign{\smallskip}\hline\noalign{\smallskip}
 0 & 40 & 6.5 & 10 & 2.44948974  & 2.44948974 & 0.00e-00 \\
 1 & 40 & 6.5 & 10 & 3.14159265 & 3.14159265 & 0.00e-00 \\
 1.5 & 30 & 2 & 4 & 3.65375388 & 3.65375374 & 1.40e-07 \\
 2 & 20 & 4 & 6 & 4.35287462 & 4.35287460 & 2.00e-08 \\
 2.5 & 20 & 4 & 6 & 5.35527531 & 5.35527546 & 1.50e-07 \\
 3 & 20 & 6.5 & 8 & 6.89684855 & 6.89684862 & 7.00e-08 \\
 4 & 20 & 14 & 16 & 14.9715808 & 14.9715463 & 3.45e-05 \\
\noalign{\smallskip}\hline
\end{tabular}
\end{table}
\begin{table}[htbp]
\caption{Comparison of $y(x)$ of standard Lane-Emden equation between present method and exact values given by Horedt\citep{Horedt}, for $m=2$.}
\label{Table. 3.}       
\centering \begin{tabular}{llll}
\hline\noalign{\smallskip}
 x & ICSRBF & Horedt & Error    \\
\noalign{\smallskip}\hline\noalign{\smallskip}
 0.00 & 1.0000000 & 1.0000000 & 0.00e+00 \\
 0.1 & 0.9983350 & 0.9983349 & 0.00e+00 \\
 0.5 & 0.9593527 & 0.9593527 & 0.00e+00 \\
 1.0 & 0.8486541 & 0.8486541 & 0.00e+00 \\
 3.0 & 0.2418240 & 0.2418241 & 1.00e-07 \\
 4.0 & 0.0488398 & 0.0488401 & 3.00e-07 \\
 4.3 & 0.0068107 & 0.0068109 & 2.00e-07 \\
 4.35 & 0.0003660 & 0.0003660 & 0.00e+00 \\
\noalign{\smallskip}\hline
\end{tabular}
\end{table}
\begin{table}[htbp]
\caption{Comparison of $y(x)$ of standard Lane-Emden equation between present method and exact values given by Horedt\citep{Horedt}, for $m=3$.}
\label{Table. 4.}       
\centering \begin{tabular}{llll}
\hline\noalign{\smallskip}
 x & ICSRBF & Horedt & Error  \\
\noalign{\smallskip}\hline\noalign{\smallskip}
 0.00 & 1.0000000 & 1.0000000 & 0.00e+00 \\
 0.1 & 0.9983358 & 0.9983358 & 0.00e+00 \\
 0.5 & 0.9598391 & 0.9598391 & 0.00e+00 \\
 1.0 & 0.8550575 & 0.8550576 & 1.00e-07 \\
 5.0 & 0.1108198 & 0.1108198 & 0.00e+00 \\
 6.0 & 0.0437379 & 0.0437380 & 1.00e-07 \\
 6.8 & 0.0041677 & 0.0041678 & 1.00e-07 \\
\noalign{\smallskip}\hline
\end{tabular}
\end{table}
\begin{table}[htbp]
\caption{Comparison of $y(x)$ of standard Lane-Emden equation between present method and exact values given by Horedt\citep{Horedt}, for $m=4$.}
\label{Table. 5.}       
\centering \begin{tabular}{llll}
\hline\noalign{\smallskip}
 x & ICSRBF  & Horedt & Error   \\
\noalign{\smallskip}\hline\noalign{\smallskip}
 0.00 & 1.0000000 & 1.0000000 & 0.00e+00 \\
 0.1 & 0.9983401 & 0.9983367 & 3.40e-06 \\
 0.2 & 0.9933921 & 0.9933862 & 5.90e-06 \\
 0.5 & 0.9603107 & 0.9603109 & 2.00e-07 \\
 1.0 & 0.8608195 & 0.8608138 & 5.70e-06 \\
 5.0 & 0.2359598 & 0.2359227 & 3.71e-06 \\
 10 & 0.05965343 & 0.0596727 & 1.92e-05 \\
 14 & 0.0083590 & 0.0083305 & 2.85e-05 \\
 14.9 & 0.00058404 & 0.0005764 & 7.64e-06 \\
\noalign{\smallskip}\hline
\end{tabular}
\end{table}
\begin{table}[htbp]
\caption{$\|Res\|_2$ of the standard Lane-Emden equations for $m=1.5, 2, 2.5, 3, 4$ respectively.}
\label{Table. 6.}       
\centering  \begin{tabular}{llllll}
\hline\noalign{\smallskip}
n & m=1.5 & m=2 & m=2.5 & m=3 & m=4  \\
\noalign{\smallskip}\hline\noalign{\smallskip}
 5 & 5.516e-03 & 1.434e-02 & 1.338e-03 & 9.996e-02 & 1.869e-00 \\
 10 & 7.250e-05 & 1.646e-04 & 7.148e-06 & 1.851e-04 & 6.436e-03 \\
 15 & 1.050e-06 & 1.379e-05 & 5.919e-09 & 1.699e-08 & 5.362e-03 \\
 20 & 1.731e-08 & 1.869e-06 & 1.268e-12 &  9.786e-10 & 3.049e-05 \\
\noalign{\smallskip}\hline
\end{tabular}
\end{table}

\begin{table}[htbp]
\caption{maximum error of the standard Lane-Emden equations for $m=0, 1, 5$ respectively.}
\label{Table. 7.}       
\centering \begin{tabular}{llll}
\hline\noalign{\smallskip}
 n & m=0 & m=1 & m=5    \\
\noalign{\smallskip}\hline\noalign{\smallskip}
 5 & 2.70199e-01 & 7.87765e-02 & 3.54822e-01 \\
 10 & 6.90811e-03 & 2.30291e-03 & 6.61006e-03 \\
 15 & 1.15061e-03 & 2.87683e-04 & 9.08909e-04 \\
 20 & 2.41206e-04 & 4.78916e-05 & 1.71788e-05 \\
\noalign{\smallskip}\hline
\end{tabular}
\end{table}

\begin{figure}[htbp!]
\center\includegraphics[scale=0.7]{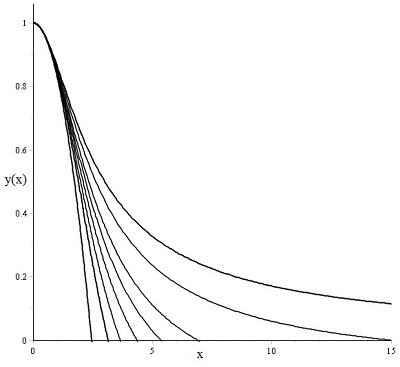}
\caption{Graph of standard Lane-Emden equation for m=0, 1, 1.5, 2, 2.5, 3, 4 and 5.}
\label{Fig. 1.}
\end{figure}

For  $p(x)=1,~q(y(x))=y^m(x),~\alpha=2,~A=1,~h(x)=0,$ and $B=0$, Eq. (\ref{lane-eq}) is the standard Lane-Emden eqution that was used to model the thermal behaviour of a spherical cloud of gas acting under the mutual attracting of its molecules and subject to the classical laws of thermodynamics\citep{Shawagfeh,Davis}.
\begin{equation}
y''(x)+\frac{2}{x}y'(x)+y^m(x)=0,   x\geq 0,
\label{Gen-LE}
\end{equation}
subject to the boundary conditions\\
\begin{equation}
y(0)=1,~~~y'(0)=0,\nonumber
\end{equation}
where $m\geq  0$ is constant. Substituting $m=0,~1,~5$ into Eq. (\ref{Gen-LE}) leads to the exact solution
\begin{eqnarray*}
&&y(x)=1-\frac{1}{3!}x^2, \\ 
&&y(x)=\frac{\sin(x)}{x}, \\  
&&y(x)=\Bigg(1+\frac{x^2}{3}\Bigg)^{-\frac{1}{2}}, \\ 
\end{eqnarray*}
respectively.\\
In other cases there aren't any analytic exact solutions. Therefore, we apply ICSRBF method to solve the standard Lane-Emden Eq. (\ref{Gen-LE}) for $m=0,~ 1,~ 1.5,~ 2,~ 2.5,~ 3,~ 4$ and $5$. To this way now we can construct the residual functions as follows:
\begin{equation}
Res(x)=xy^{(2)}_n(x)+2 (y^{(1)}_n(x))+x(y_n(x))^m,
\end{equation}
As said before, to obtain the coefficient $\xi_is,~Res(x)$ is equalized to zero at $N$ points by (\ref{formula}) :
\begin{equation}
Res(x_j)=0,~~~j=0, 1, 2, ..., N.\nonumber
\end{equation}
By solving this set of equations, we can find the approximating function $y_n(x)$.\\
Table \ref{Table. 2.} shows the comparison of the first zeros of  standard Lane-Emden equations, from the present method and exact given by Horedt\citep{Horedt} for $m=0,~ 1,~ 1.5,~ 2,~ 2.5,~ 3$ and 4, respectively.\\
Table \ref{Table. 3.}, \ref{Table. 4.}, \ref{Table. 5.} shows the approximation of $y(x)$ for the standard Lane-Emden equation for $m=2,~ 3,~4$ respectively obtained by the method proposed in this paper and those obtained by Horedt\citep{Horedt}.\\
Table \ref{Table. 6.} represents the $\|Res\|_2$ by the present method for $m=1.5,~ 2,~2.5,~ 3,~ 4$ for several points and Table \ref{Table. 7.} represent the maximum error by the present method for $m=0, ~1,~ 5$ for several points to how that the new method has appropriate convergence rate. The result graph of the standard Lane-Emden equation for $m=0, ~1, ~1.5, ~2, ~2.5, ~3, ~4$ and 5 is shown in Fig. \ref{Fig. 1.}.\\

\subsection{Example 2 (The isothermal gas sphere equation)}
\label{sssec52}
\begin{figure}
\center\includegraphics[scale=0.7]{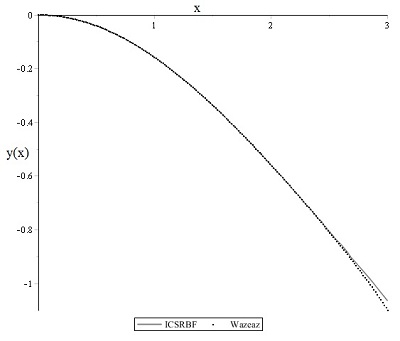}
\caption{\scriptsize Graph of Isothermal gas sphere equation in comparison with Wazwaz solution\citep{Wazwaz1}}
\label{Fig. 2.}
\end{figure}

\begin{figure}
\center\includegraphics[scale=0.7]{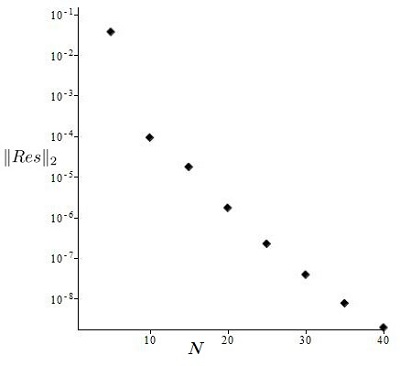}
\caption{\scriptsize The $\|Res\|_2$  of Isothermal gas sphere equation in for N=5, 10, 15, 20, 25, 30, 35 and 40}
\label{Fig. 3.}
\end{figure}

\begin{table}[htbp]
\caption{Comparison of $y(x)$, between present method and series solution given by Wazwaz\citep{Wazwaz1} and HFC method\citep{Parand7} for isothermal gas sphere.}
\label{Table. 8.}       
\centering \begin{tabular}{llll}
\hline\noalign{\smallskip}
 x  & ICSRBF & HFC & Adomain  \\
\noalign{\smallskip}\hline\noalign{\smallskip}
 0.00 & 0.0000000000 & 0.0000000000 & 0.0000000000 \\
 0.1 & -0.0016658338 & -0.0016664188 & -0.0016658339 \\
 0.2 & -0.006533671 & -0.0066539713 & -0.0066533671 \\
 0.5 & -0.0411539573 & -0.0411545150 & -0.0411539568 \\
 1.0 & -0.1588276775 & -0.1588281737 & -0.1588273537 \\
 1.5 & -0.3380194248 & -0.3380198308 & -0.3380131103 \\
 2.0 & -0.5598230044 & -0.5598233120 & -0.5599626601 \\
 2.5 & -0.8063408707 & -0.8063410846 & -0.8100196713 \\
\noalign{\smallskip}\hline
\end{tabular}
\end{table}
For $p(x)=1,~q(y(x))=e^{y(x)},~h(x)=0,~\alpha=2,~A=0$ and $B=0$, Eq.(\ref{lane-eq}) is the isothermal gas sphere equation
\begin{equation}
\label{Eq17}
y''(x)+\frac{2}{x}y'(x)+e^{y(x)}=0,   x\geq 0,
\end{equation}\\
subject to the boundary conditions\\
\begin{equation}
y(0)=0,~~~y'(0)=0,\nonumber\\
\end{equation}
 Wazwaz\citep{Wazwaz1} by using ADM has produced a series solution as follow:
\begin{equation}
\label{Eq18}
y(x)\simeq -\frac{1}{6}x^2+\frac{2}{5.4!}x^4-\frac{8}{21.6!}x^6+\frac{122}{81.8!}x^8-\frac{61.67}{495.10!}x^{10}.
\end{equation}
We applied the ICSRBF method to solve this equation (\ref{Eq17}). We construct the residual function as follows:\\
\begin{equation}
Res(x)=xy^{(2)}_n(x)+2 (y^{(1)}_n(x))+xe^{(y_n(x))},
\end{equation}
To obtain the coefficients $\xi_is,~Res(x)$ is equalized to zero at $N$ point by (\ref{formula}) with $\gamma=1.7$:\\
\begin{equation}
Res(x_j)=0,~~~j=0,1,2,...,N.\nonumber
\end{equation}
By solving this set of equations, we have the approximating function $y(x)$. Table \ref{Table. 8.} shows the comparison of $y(x)$ obtained by method proposed in this paper with $(N=40,~r_\omega=6.5)$, and those obtained by Wazwaz\citep{Wazwaz1} (\ref{Eq18}) and HFC method by Parand et al\citep{Parand7}.\\
In order to compare the present method with those obtained by Wazwaz\citep{Wazwaz1} the resulting graph of Eq.(\ref{Eq16}) is shown in Fig. \ref{Fig. 2.}. The graph of residual for this equation for $N= 5, ~10, ~15, ~20, ~25, ~30, ~35$ and 40 is shown in Fig. \ref{Fig. 3.}.

\subsection{Example 3 (The White Dwarf equation)}
\label{sssec53}
\begin{table}
\caption{Comparison of $y(x)$, between present method and series solution given by Singh\citep{Singh} and Haar method\citep{Kaur} for White Dwarf equation.}
\label{Table. 9.}       
\centering \begin{tabular}{lllll}
\hline\noalign{\smallskip}
 x &$\sigma$& ICSRBF & MHAM & Haar    \\
\noalign{\smallskip}\hline\noalign{\smallskip}
 0.0001 & 0.1 & 0.999999 & 0.999999 & 1 \\
 0.01 &   & 0.999985 & 0.999985 & 0.999986 \\
 0.1 &   & 0.998579 & 0.998578 & 0.998581 \\
 0.2 &   & 0.994340 & 0.994340 & 0.994379 \\
 0.4 &   & 0.977738 & 0.977738 & 0.978345 \\
 0.6 &   & 0.951270 & 0.951263 & 0.953273 \\
 0.7 &   & 0.934823 & 0.934801 & 0.936497 \\
 0.9 &   & 0.896673 & 0.896522 & 0.895488 \\
\noalign{\smallskip}
 0.0001 & 0.2 & 0.999999 & 0.999999 & 1 \\
 0.01 &   & 0.999988 & 0.999988 & 0.999988 \\
 0.1 &   & 0.998809 & 0.998809 & 0.998580 \\
 0.2 &   & 0.995255 & 0.995251 & 0.995296 \\
 0.4 &   & 0.981320 & 0.981320 & 0.981931 \\
 0.6 &   & 0.959049 & 0.959045 & 0.961128 \\
 0.7 &   & 0.945179 & 0.945165 & 0.947116 \\
 0.9 &   & 0.912913 & 0.912812 & 0.912833 \\
\noalign{\smallskip}
 0.0001 & 0.3 & 0.999999 & 0.999999 & 1 \\
 0.01 &   & 0.999990 & 0.999990 & 0.999999 \\
 0.1 &   & 0.999025 & 0.999025 & 0.999028 \\
 0.2 &   & 0.996115 & 0.996115 & 0.996154 \\
 0.4 &   & 0.984690 & 0.984690 & 0.985295 \\
 0.6 &   & 0.966384 & 0.966381 & 0.968464 \\
 0.7 &   & 0.954953 & 0.954944 & 0.956942 \\
 0.9 &   & 0.928280 & 0.928216 & 0.928599 \\
\noalign{\smallskip}\hline
\end{tabular}
\end{table}

\begin{figure}
\center\includegraphics[scale=0.7]{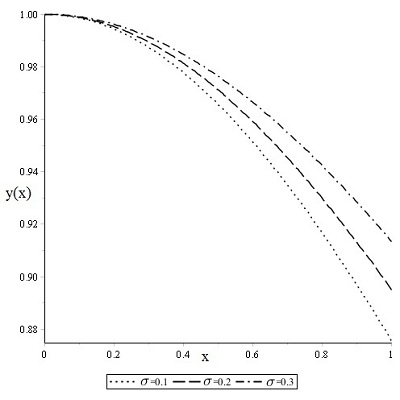}
\caption{\scriptsize Graph of White Dwarf equation for $\sigma$=0.1, 0.2 and 0.3}
\label{Fig. 4.}
\end{figure}

For $p(x)=1,~q(y(x))=(y^2-\sigma)^{\frac{3}{2}},~h(x)=0,~\alpha=2,~A=1$ and $B=0$, Eq.(\ref{lane-eq})\\
will be one of the white Dwarf
 equation that is absorbing to solve
\begin{equation}
\label{Eq15}
y''(x)+\frac{2}{x}y'(x)+(y(x)^2-\sigma)^{\frac{3}{2}}=0,   x\geq 0,
\end{equation}\\
subject to the boundry conditions\\
\begin{equation}
y(0)=1,~~~y'(0)=0,\nonumber
\end{equation}
which was introduced by Chandraskhar \citep{Chandrasekhar} in his study of the gravitational potential of the degenerate White Dwarf stars.
A series solution obtained by Singh\citep{Singh} by using MHAM as follows:
\begin{equation}
\label{Eq16}
y(x)\simeq 1-\frac{1}{6}\varpi^3x^2+\frac{1}{40}\varpi^4x^4-\frac{1}{7!}\varpi^5[5\varpi^2+14]x^6.
\end{equation}
where $\varpi=\sqrt{1-\sigma}$.\\
We applied the ICSRBF method to solve this equation (\ref{Eq15}). We construct the residual function as follows:\\
\begin{equation}
Res(x)=xy^{(2)}_n(x)+2 (y^{(1)}_n(x))+x({(y_n(x))}^2-\sigma^2)^\frac{3}{2},
\end{equation}
To obtain the coefficients $\xi_is,~Res(x)$ is equalized to zero at $N$ point by (\ref{formula}) with $\gamma=1.5$:\\
\begin{equation}
Res(x_j)=0,~~~j=0,1,2,...,N.\nonumber
\end{equation}
By solving this set of equations,we have the approximating function $y(x)$. Table \ref{Table. 9.} shows the comparison of $y(x)$ obtained by method proposed in this paper with $(N=20,~r_\omega=0.5)$, and those obtained by singh\citep{Singh} (\ref{Eq16}) and Haar method by Kaur\citep{Kaur}.\\
 The result graph of the White Dwarf equation for $\sigma$= 0.1, 0.2 and 0.3  is shown in Fig. \ref{Fig. 4.}.\\

\subsection{Example 4}
\label{sssec54}
\begin{table}
\caption{Comparison of $y(x)$, between present method and series solution given by Wazwaz\citep{Wazwaz1} and HFC method\citep{Parand7} for Example 4.}
\label{Table. 10.}       
\centering \begin{tabular}{llll}
\hline\noalign{\smallskip}
 x & ICSRBF & HFC & Adomain   \\
\noalign{\smallskip}\hline\noalign{\smallskip}
 0.00 & 1.0000000000 & 1.0000000000 & 1.0000000000 \\
 0.1 & 0.9980430038 & 0.9981138095 & 0.9980428414 \\
 0.2 & 0.9921896287 & 0.9922758837 & 0.9921894348 \\
 0.5 & 0.9519612468 & 0.9520376245 & 0.9519611019 \\
 1.0 & 0.8182430031 & 0.8183047481 & 0.8182516669 \\
 1.5 & 0.6254386159 & 0.6254886192 & 0.6258916077 \\
 2.0 & 0.4066218732 & 0.4066479695 & 0.4136691039 \\
\noalign{\smallskip}\hline
\end{tabular}
\end{table}
For $p(x)=1,~q(y(x))=\sinh(y(x)),~h(x)=0,~\alpha=2,~A=1$ and $B=0$, Eq.(\ref{lane-eq}) will be one of the Lane-Emden type equations that is absorbing to solve
\begin{equation}
\label{Eq13}
y''(x)+\frac{2}{x}y'(x)+\sinh(y)=0,   x\geq 0,
\end{equation}\\
subject to the boundry conditions\\
\begin{equation}
y(0)=1,~~~y'(0)=0,\nonumber
\end{equation}

\begin{figure}
\center\includegraphics[scale=0.7]{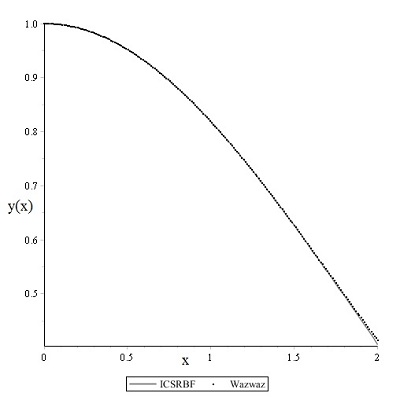}
\caption{\scriptsize Graph of equation of example 4 in comparing the present method and Wazwaz solution \citep{Wazwaz1}.}
\label{Fig. 5.}
\end{figure}

\begin{figure}
\center\includegraphics[scale=0.7]{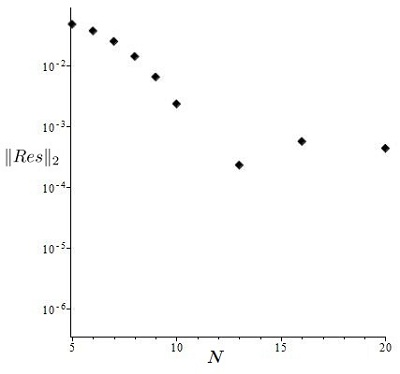}
\caption{\scriptsize The $\|Res\|_2$   of example 4  for $N= 5,~ 6,~ 7,~ 8,~ 9,~ 10,~ 13,~ 16$ and 20 }
\label{Fig. 6.}
\end{figure}

A series solution obtained by Wazwaz\citep{Wazwaz1} by using ADM is:
\begin{eqnarray}\footnotesize
\label{Eq14}
&&y(x)\simeq 1-\frac{(e^2-1)x^2}{12e}+\frac{1}{480}\frac{(e^4-1)x^4}{e^2}-\nonumber\\
&&\frac{1}{30240}\frac{(2e^6+3e^2-3e^4-2)x^6}{e^3}+\nonumber\\
&&\frac{1}{26127360}\frac{(16e^8-104e^6+104e^2-61)x^8}{e^4}.
\end{eqnarray}
We applied the ICSRBF method to solve this equation (\ref{Eq13}). We construct the residual function as follows:\\
\begin{equation}
Res(x)=xy^{(2)}_n(x)+2(y^{(1)}_n(x))+x\sinh({y_n(x)}),
\end{equation}
To obtain the coefficients $\xi_is,~Res(x)$ is equalized to zero at $N$ point by (\ref{formula}) with $\gamma=1.7$:\\
\begin{equation}
Res(x_j)=0,~~~j=0,1,2,...,N.\nonumber
\end{equation}
By solving this set of equations,we have the approximating function $y(x)$. Table \ref{Table. 10.} shows the comparison of $y(x)$ obtained by method proposed in this paper with $(N=20,~r_\omega=1)$,and those obtained by Wazwaz\citep{Wazwaz1} (\ref{Eq14}) and HFC method by Parand et al.\citep{Parand7}.\\
In order to compare the present method with those obtained by Wazwaz\citep{Wazwaz1} the resulting graph of Eq.(\ref{Eq16}) is shown in Fig. \ref{Fig. 5.}. The graph of $\|Res\|_2$ for this equation for $N= 5, ~6, ~7,~ 8,~ 9,~ 10,~ 13,~ 16$ and 20 is shown in Fig. \ref{Fig. 6.}.\\

\subsection{Example 5}
\label{sssec55}
\begin{figure}
\center\includegraphics[scale=0.7]{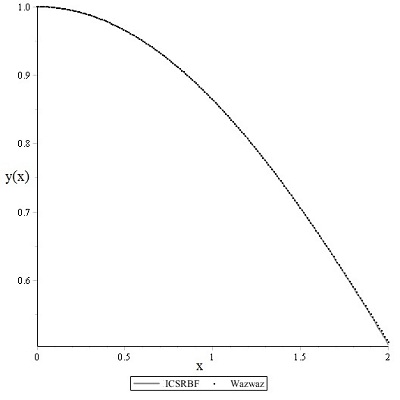}
\caption{\scriptsize Graph of equation of example 5 in comparing the present method and Wazwaz solution \citep{Wazwaz1}.}
\label{Fig. 7.}
\end{figure}

\begin{figure}
\center\includegraphics[scale=0.7]{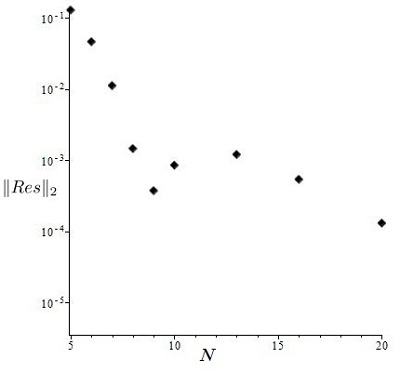}
\caption{\scriptsize The $\|Res\|_2$  graph of example 5  for $N= 5,~ 6,~ 7,~ 8,~ 9,~ 10,~ 13,~ 16$ and 20 }
\label{Fig. 8.}
\end{figure}

\begin{table}
\caption{Comparison of $y(x)$, between present method and series solution given by Wazwaz\citep{Wazwaz1} and HFC method\citep{Parand7} for Example 5.}
\label{Table. 11.}       
\centering \begin{tabular}{llll}
\hline\noalign{\smallskip}
 x & ICSRBF & HFC & Adomian   \\
\noalign{\smallskip}\hline\noalign{\smallskip}
 0.0 & 1.0000000000 & 1.0000000000 & 1.0000000000 \\
 0.1 & 0.9985979436 & 0.9986051425 & 0.9985979358 \\
 0.2 & 0.9943962892 & 0.9944062706 & 0.9943962733 \\
 0.5 & 0.9651778048 & 0.9651881683 & 0.9651777886 \\
 1.0 & 0.8636811302 & 0.8636881301 & 0.8636811027 \\
 1.5 & 0.7050451897 & 0.7050524103 & 0.7050419247 \\
 2.0 & 0.5064632371 & 0.5064687568 & 0.5063720330 \\
\noalign{\smallskip}\hline
\end{tabular}
\end{table}
For $p(x)=1,~q(y(x))=\sin(y(x)),~h(x)=0,~\alpha=2,~A=1$ and $B=0$, Eq.(\ref{lane-eq}) will be one of the Lane-Emden type equations that is absorbing to solve
\begin{equation}
\label{Eq11}
y''(x)+\frac{2}{x}y'(x)+\sin(y)=0,   x\geq 0,
\end{equation}\\
subject to the boundry conditions\\
\begin{equation}
y(0)=1,~~~ y'(0)=0,\nonumber
\end{equation}
A series solution obtained by Wazwaz\citep{Wazwaz1} by using ADM is:
\begin{equation}
\label{Eq12}
y(x)\simeq 1-\frac{1}{6}k_1x^2+\frac{1}{120}k_1k_2x^4+k_1(\frac{1}{3024}k_1^2.-\frac{1}{5040}k_2^2)x^6+
k_1k_2(-\frac{113}{3265920}k^2_1+ \\  \frac{1}{362880}k^2_2)x^8+...
\end{equation}
where $k_1=\sin(1)$ and $k_2=\cos(1)$.\\
We applied the ICSRBF method to solve this equation (\ref{Eq11}). We construct the residual function as follows:\\
\begin{equation}
Res(x)=xy^{(2)}_n(x)+2(y^{(1)}_n(x))+x\sin({y_n(x)}),
\end{equation}
To obtain the coefficients $\xi_is,~Res(x)$ is equalized to zero at $N$ point by (\ref{formula}) with $\gamma=1.6$:\\
\begin{equation}
Res(x_j)=0,~~~j=0,1,2,...,N.\nonumber
\end{equation}
By solving this set of equations, we have the approximating function $y(x)$. Table \ref{Table. 11.} shows the comparison of $y(x)$ obtained by method proposed in this paper with $(N=20,~r_\omega=2)$, and those obtained by Wazwaz\citep{Wazwaz1} (\ref{Eq12}) and HFC method by Parand et al. \citep{Parand7}.\\
In order to compare the present method with those obtained by Wazwaz\citep{Wazwaz1} the resulting graph of Eq.(\ref{Eq11}) is shown in Fig. \ref{Fig. 7.}.\\ The graph of $\|Res\|_2$ for this equation for $N= 5, ~6, ~7, ~8, ~9, ~10, ~13, ~16$ and 20 is shown in Fig. \ref{Fig. 8.}.\\

\section{ Conclusion}
\label{sec:6}
 Lane-Emden equations occur in the theory of stellar structure and describe the temperature variation of a spherical gas cloud. White Dwarf equation appears in the gravitational potential of the degenerate White Dwarf stars. Lane-Emden type equations have been considered by many mathematicians as mentioned before\citep{Dehghan3}. The fundamental goal of this paper has been to construct an approximation to the solution of non-linear Lane-Emden type equation in a semi-infinite interval. CSRBFs are proposed to provide an effective but simple way to improve the convergence of the solution by the collocation method. The validity of the method is based on the assumption that it converges by increasing the number of collocation points. A comparison is made among the numerical solution of Horedt\citep{Horedt} and series solutions of Wazwaz\citep{Wazwaz1}, Singh et al.\citep{Singh}, Kaur et al.\citep{Kaur} and the current work. It has been shown that our present work provides an acceptable approach for the Lane-Emden type equations. Also it was confirmed by maximum error and $\|Res\|_2$ figures, this approach has an exponentially convergence rate. Also, in this paper we show that the ICSRBF method for solving ordinary differential equations is simple and it has high accuracy and reliable convergence. \\




\bibliographystyle{elsarticle-harv} 
\bibliography{paperen}





\end{document}